\documentclass[10pt]{amsart}

\usepackage[all, cmtip]{xy}
\usepackage{amsmath,amsthm,amssymb,amsfonts,amscd}
\usepackage{tikz-cd}
\usepackage{float} 

\usepackage{graphicx} 
 
\newtheorem{theorem}{Theorem}[section]
\newtheorem{lemma}[theorem]{Lemma}
\newtheorem{proposition}[theorem]{Proposition}
\newtheorem{corollary}[theorem]{Corollary}

\theoremstyle{definition}
\newtheorem{definition}[theorem]{Definition}

\newtheorem{example}[theorem]{Example}
\newtheorem{conjecture}[theorem]{Conjecture}

\newtheorem{question}[theorem]{Question}

\theoremstyle{remark}
\newtheorem{remark}[theorem]{Remark}

\numberwithin{equation}{section}

\begin{document}

\title[On Barycentric transformations of Fano polytopes]{On Barycentric transformations \\ of  Fano polytopes}

\author{DongSeon Hwang}
\address{Department of Mathematics, Ajou University, Suwon 16499, Republic Of Korea}
\email{dshwang@ajou.ac.kr}

\author{Yeonsu Kim}
\email{kys3129@ajou.ac.kr}

\date{November 3, 2020}



\begin{abstract}
We introduce the notion of barycentric transformation of Fano polytopes, from which we can assign a certain type to each Fano polytope.  The type can be viewed as a measure of the extent to which the given Fano polytope is close to be Kähler-Einstein. In particular, we expect that every K\"ahler-Einstein or symmetric Fano polytope is of type $B_\infty$. We verify this expectation for some low dimensional cases. We emphasize that for a Fano polytope $X$  of dimension $1$, $3$ or $5$, $X$ is K\"ahler-Einstein if and only if it is of type $B_\infty$.
\end{abstract}

\maketitle




\section {Introduction}

A \emph{Fano polytope} of dimension $n$ is a full dimensional convex lattice polytope in $\mathbb{R}^n$  such that the vertices are primitive lattice points and the origin is an interior point. Note that the class of Fano polytopes of dimension $n$ up to unimodular transformation has a one-to-one correspondence with  the class of toric Fano varieties of dimension $n$ up to isomorphism. A Fano polytope $P$ is said to be \emph{K\"ahler-Einstein} if the dual polytope of $P$ has the origin as its barycenter. By Theorem \cite[Theorem 1.2]{BB}, a Fano polytope $P$ is  K\"ahler-Einstein if and only if the associated toric Fano variety $X_P$ admits a K\"ahler-Einstein metric.

Batyrev and Selivanova introduced the notion of symmetric Fano polytopes to study K\"ahler-Einstein polytopes. A Fano polytope $P$ is said to be  \emph{symmetric} if the origin is the only lattice point fixed by every automorphisms of $P$ onto itself.

\begin{theorem}\cite[Theorem 1.1]{BS} 
Every smooth symmetric Fano polytope is  K\"ahler-Einstein. 
\end{theorem}

There had been interests on whether the converse statement holds(\cite{BS}, \cite[Remark in p.1257]{So}, \cite[Remark 4.3]{CL}, \cite[p.257]{FOS}). For smooth Fano polytopes of dimension at most $8$, an exhaustive investigation yields the following theorem.  

\begin{theorem}\cite[Proposition 2.1]{NP}\label{NP-theorem}
Let $P$ be a smooth K\"ahler-Einstein Fano polytope of dimension at most $8$, then  $P$ is not symmetric  if and only if $P$ is one of the three Fano polytopes $Q_1,Q_2$ or $Q_3$ in \cite{NP} where $Q_2$ and $Q_3$ have dimension $8$ and $Q_1$ has dimension $7$.
\end{theorem}

We remark that, for each integer $n \geq 9$, there exists a smooth  K\"ahler-Einstein Fano polytope of dimension $n$ that is not symmetric(\cite[Corollary 5.3]{N2}). See also \cite{HKim} for the discussion in the singular setting.

In this note, we shall describe a class of Fano polytopes that is expected to be very similar to the class of K\"ahler-Einstein Fano polytopes in the smooth case and have a close connection to  symmetric or K\"ahler-Einstein Fano polytopes in general.

To be more precise,   we shall introduce the notion of barycentric transformation, or B-transformation in short, of Fano polytopes. See Definition \ref{B-transformation} for the precise definition. The B-transformation of a Fano polytope does not always produce a Fano polytope as we see in Example \ref{badbehavior}. By using this phenomenon of B-transformation, we define a \emph{type} of Fano polytopes.  Naively speaking, a Fano polytope is \emph{of type $B_k$} if it is possible to take  the B-transformation $k$ times, and is \emph{of strict type $B_k$} if we can take the B-transformation at most $k$ times. See Definition \ref{type}  for the precise definition.
  
\begin{theorem}\label{main-smooth}
Let $P$ be a smooth  Fano polytope of dimension $n$. 
    \begin{enumerate} 
    \item Assume that $n$ is odd and $n \leq 5$.  Then, $P$ is K\"ahler-Einstein if and only if $P$ is  of type $B_\infty$. 
    \item Assume that $n$ is even and $n \leq 4$.  If $P$ is K\"ahler-Einstein, then $P$ is of type $B_\infty$ except possibly for one case with GRDB ID 54 (\cite{GRDB} based on \cite{O}).
    \item If $P$ is a K\"ahler-Einstein Fano polytope of dimension $7$ or $8$ that is not symmetric, i.e.,  $P=Q_1, Q_2,Q_3$ in \cite{NP}, then $P$ is of  type $B_\infty$. 
    \end{enumerate}
\end{theorem}

Moreover, the K\"ahler-Einstein property is preserved by the B-transformation in low dimension.

\begin{theorem}\label{BpreservesKE}
Let $P$ be a smooth Fano polytope of dimension at most $6$, or  $P=Q_1, Q_2,Q_3$ in \cite{NP}. If $P$ is K\"ahler-Einstein, then so is   $B(P)$.
\end{theorem}

Based on these results, we boldly pose the  following conjecture.

\begin{conjecture}\label{mainconj}
Let $P$ be a smooth Fano polytope. 
    \begin{enumerate}
        \item If $P$ has odd dimension, then $P$ is K\"ahler-Einstein if and only if $P$ is of type $B_\infty$.
        \item If $P$ is K\"ahler-Einstein, then $B(P)$ is also a K\"ahler-Einstein Fano polytope. In particular,  $P$ is of type $B_\infty$. 
        \item If $P$ is symmetric, then $B(P)$ is a symmetric and K\"ahler-Einstein Fano polytope. In particular, $P$ is of type $B_\infty$.
    \end{enumerate}
\end{conjecture}
See Question \ref{pseudoperiodic} for the discussion on the converse statement of Conjecture \ref{mainconj} (2) in the even dimensional case. 

One can ask whether  Conjecture \ref{mainconj} holds true for the singular varieties.   Contrary to the smooth case, not every symmetric Fano polytope is K\"ahler-Einstein (\cite{HKim}) even for surfaces. But we still have the following result.

\begin{theorem}\label{main-singular}
Let $P$ be a Fano polygon.
\begin{enumerate}
    \item If $P$ is K\"ahler-Einstein, then $P$ is of type  $B_1$. 
    If, moreover, $P$ is a triangle, then $B(P)$ is  a K\"ahler-Einstein Fano triangle, hence $P$ is of type  $B_\infty$.
    \item If $P$ is symmetric, then $B(P)$ is a symmetric and K\"ahler-Einstein Fano polygon, hence $P$ is of type $B_\infty$. 
\end{enumerate}
\end{theorem}

Every known example of a K\"ahler-Einstein Fano polygon that is not symmetric is a triangle. For this reason, the following was suggested.

\begin{conjecture}\cite[Conjecture 1.6]{HKim}\label{triangle}
Let $P$ be a K\"ahler-Einstein Fano polygon. If $P$ is not symmetric, then it is a triangle.
\end{conjecture}
Assuming Conjecture \ref{triangle}, Theorem \ref{main-singular} proves Conjecture \ref{mainconj} for surfaces. Theorem \ref{index17} shows that Conjecture \ref{triangle}, hence Conjecture \ref{mainconj}, holds true for Fano polygons of index at most $17$ without any assumption.
 
\begin{remark}
Not every Fano polytope of type $B_\infty$ becomes K\"ahler-Einstein after a suitable application of B-transformations. For example, consider a Fano polygon $P$ with GRDB ID 13118 (\cite{GRDB} based on \cite{KKN}). Then, $P$ is $2$-periodic, i.e., $B^2(P) = P$. But both $P$ and $B(P)$ are not K\"ahler-Einstein. Note that $P$ is of index $3$, which is the smallest possible index with this property. 
\end{remark}

Section 2 provides rigorous definitions and examples  of barycentric transformations. Main theorems are proven in Section 3. Section 4 presents some interesting observation on the orbits of a Fano polygon under B-transformations and discussion on Fano polygons with zero barycenter.

\section{B-transformation: definition and examples} 

\subsection{Notation} 
For a lattice point $v = (x, y)$, we define the \emph{primitive index} $I(v)$ of $v$ by $I(v)=gcd(x,y)$. A lattice point is called \emph{primitive} if its primitive index is one. 
The \emph{order} of the two dimensional cone $C$ spanned by two lattice points $v_i = (x_i, y_i)$ and $v_{i+1} = (x_{i+1}, y_{i+1})$, denoted by $ord(v_i,v_{i+1})$, is defined by 
$$ord(v_i,v_{i+1}):=det\left( \begin{array}{ccr} x_i & x_{i+1} \\
y_i & y_{i+1} \\ \end{array} \right)=x_iy_{i+1}-y_ix_{i+1}.$$

\subsection{B-transformations} 
To understand more about K\"ahler-Einstein Fano polytopes, we introduce the new notion called a barycentric  transformation.

\begin{definition} \label{B-transformation} 
For a Fano polytope $P$, we can always obtain a lattice polytope  by taking the convex hull of the  barycenters of all maximal dimensional cones of $P$. See (\cite[Exercise 11.1.10]{CLS}) for the definition of the barycenter of a cone. This association is called a \emph{barycentric transformation} of $P$ or, in short, a \emph{B-transformation} of $P$. 
\end{definition}

Note that the B-transformation is uniquely determined up to unimodular transformation.

We can easily compute the B-transformation of a Fano polygon using the following lemma, which immediately follows from definition.

\begin{lemma} 
Let $P$ be a Fano polygon with vertices $v_1,\dots,v_n$ written in counterclockwise order. 
Then,  $$B(P)=conv\Big\{\frac{v_1+v_2}{I(v_1+v_2)}, \dots,\frac{v_n+v_1}{I(v_n+v_1)} \Big\}.$$
\end{lemma}

\begin{example} \label{badbehavior}
For a Fano polytope $P$, $B(P)$ is not a Fano polytope in general.  
    \begin{enumerate} 
        \item Take $P=conv\{(2,-1),(0,1),(-1,0)\}$.         
        Then, $$B(P)=conv\{(1,0),(-1,1),(1,-1)\}$$  
        is not a Fano polytope since it does not contain the origin as an interior point. 
        \item Take $P=conv\{(1, -2), (0, 1), (-1, -2)\}.$
        Then, $$B(P)=conv\{(1,-1),(-1,-1)\}$$ is not a Fano polytope. Note that the dimension is decreased in this case. 
    \end{enumerate}
\end{example}

\begin{definition} \label{type}
Let $P$ be a Fano polytope of dimension $d$.  
\begin{enumerate} 
    \item $P$ is said to be \emph{of type $B_m$} if $B^m(P)$ is a Fano polytope of dimension $d$.
    \item $P$ is said to be \emph{of strict type $B_m$} if it is of type $B_m$ but not of type $B_{m+1}$. 
    \item $P$ is said to be \emph{of type $B_\infty$} if it is of type $B_m$ for every positive integer $m$. 
\end{enumerate}
By convention, every Fano polytope $P$ is of type $B_0$. Note that if a Fano polytope $P$ is type of $B_s$ then it is of type $B_t$ for every $t < s$.
\end{definition}

We also introduce the following useful notion.

\begin{definition} Let $P$ be a Fano polytope.
    \begin{enumerate}
    \item $P$ is said to be \emph{B-invariant} if $B(P) = P$.  
        \item  $P$ is called $k$-\emph{periodic} (or just \emph{periodic}) if there exists an integer $t$ such that $B^{k+t}(P) = B^t(P)$ for some integer $k \geq 1$. 
        \item $P$ is called \emph{pseudo-periodic} if the number of vertices of $B^k(P)$ is invariant under the B-transformation for some positive integer $k$.
    \end{enumerate} 
\end{definition}
Note that a periodic Fano polytope is of type $B_\infty$. 

\begin{example}\label{cube}
Let $P$ be a $d$-dimensional Fano cube. For example, if $d=2$, 
$$P = conv \{ e_1+e_2, -e_1+e_2, -e_1-e_2, e_1-e_2 \}$$
where $e_1$ and $e_2$ are standard bases of $\mathbb{R}^2$. Then, $P$ is of type $B_\infty$.  
Indeed, one can easily see that $B(P)=conv\{e_1,\dots,e_d,-e_1,\dots,-e_d\}$ is a  $d$-dimensional Fano bipyramid and  $B^2(P)=P$,i.e., $P$ is $2$-periodic. In this case, both $P$ and $B(P)$ are  K\"ahler-Einstein. Note that  $P$ is singular and $B(P)$ is smooth, from which we see that   the B-transformation does not preserve the smoothness.
\end{example}

\begin{example}  
Let $P=conv\{(0,1),(3,-2),(-4,1)\}$. 
It is easy to compute that  $B(P) = conv\{(3,-1),(-1,-1),(-2,1)\} \ \text{and}$
$B^2(P) = conv\{(-1,0),(1,-1),(1,0)\},$ 
from which we see that   $P$ is of strict type $B_1$.
\end{example}

\begin{example}\label{toKE} 
Let $P=conv\{(-25,-12),(-5,-6),(25,14)\}.$ 
One can compute that $B(P)=conv\{(-5,-3),(5,2),(0,1)\}$ and 
 $$B^{2}(P)=conv\{(0,-1),(5,3),(-5,-2)\} = -B(P),$$ thus $P$ is $2$-periodic, and hence it is of type $B_\infty$. Now, it is easy to see that $P$ is not  K\"ahler-Einstein  but so is $B(P)$. Note also that $B^n(P)$ is not symmetric for every integer $n \geq 0$.  
\end{example}

\begin{lemma}\label{monotone}
For a Fano polygon, the number of vertices is not increasing under  B-transformations. 
\end{lemma}

\begin{proof}
It immediately follows from the fact that, in dimension two,   the number of maximal dimensional cones  is equal to the number of vertices.  
\end{proof}

\begin{example}  
The number of vertices is not invariant under B-transformation in general. Let $$P=conv\{(3,-1),(3,1),(1,2),(-3,1),(-3,-1),(-1,-2)\}.$$ Then, $P$ is both symmetric and K\"ahler-Einstein but $$B(P)=conv\{(4,3),(-2,3),(-4,-3),(2,-3)\}.$$ 
By using the argument in the proof of Lemma \ref{1-refl}, it is easy to see  that  for a K\"ahler-Einstein Fano polygon with at most $4$ vertices, the number of vertices is invariant under  $B$-transformation.
We do not know any example of such a K\"ahler-Einstein Fano polygon with $5$ vertices. 
\end{example}

\begin{remark}
Lemma \ref{monotone} does not hold for higher dimensional Fano polytopes. See Example \ref{cube}. 
\end{remark}

\subsection{Condition to be of type $B_1$}
We provide a sufficient condition for a given Fano polygon $P$ to be of type $B_1$.

\begin{lemma}\label{contain0}
Let $P$ be a  convex lattice polygon with primitive vertices $v_1,\dots,v_n$ written in counterclockwise order. Then, 
$P$ contain the origin as an interior point  if and only if $ord(v_i,v_{i+1})>0$ for all $i=1, \dots, n$.
\end{lemma}

\begin{proof}
If $P$ contain the origin as an interior point, then $P$ is a Fano polygon. 
By the choice of the orientation of vertices, we have  $ord(v_i,v_{i+1})>0$ for every $i$.
Conversely, suppose that $ord(v_i,v_{i+1})>0$ for all $i=1, \dots, n$. Assume that $P$ does not contain the origin as an interior point. Take a face $F$ of $P$ that is closest to the origin. Then, $F$ has two vertices $v_{i-1}$ and $v_{i}$ of $P$. Since $v_i$ is a primitive point, we can map $v_i$ to $(0,1)$ by an orientation-preserving unimodular transformation.
Then, $v_{i-1}=(a,-b)$ and $v_{i+1}=(-c,d)$ for some positive integers $a,b,c$ and $d$. Since the origin is not an interior point of  $P$,  $v_{i+1}$ is above the line generated by $v_{i-1}$ and $v_{i}$. Indeed, if otherwise, it is easy to see that there exists another face of $P$ that is closer to the origin. Now, the three vertices $v_{i-1},v_{i},v_{i+1}$ are  in clockwise order, which is a contradiction. 
\end{proof}

Let $P$ be a Fano polygon with  vertices $v_1,\dots,v_n$ written in counterclockwise order. Then, for each $i$ with $1 \leq i \leq n$, we consider the following function 
$$g(i):=ord(v_{i-1},v_{i})+ord(v_{i},v_{i+1})-ord(v_{i+1},v_{i-1}).$$

\begin{proposition}\label{cond_B_1}
Let $P$ be a Fano polygon with  vertices $v_1,\dots,v_n$ written in counterclockwise order. 
If $g(i)>0$ for all $i =1, \dots ,n$, then $P$ is of type $B_1$.
\end{proposition}

\begin{proof} 
It is enough to show that $B(P)$ contains the origin as an interior point. 
Since $g(i)>0$ for all $i=1, \dots, n$,   it is easy to see that  $$g(i)=ord(v_{i-1},v_{i})+ord(v_{i},v_{i+1})-ord(v_{i+1},v_{i-1})=ord(v_{i-1}+v_{i},v_{i}+v_{i+1})>0$$  
Note that $ord(v_{i-1}+v_{i},v_{i}+v_{i+1})>0$ if and only if $ord(\frac{v_{i-1}+v_{i}}{I_p(v_{i-1}+v_i)},\frac{v_{i}+v_{i+1}}{I_p(v_i+v_{i+1})})>0$. Then, the origin is an interior point of $B(P)$ by Lemma \ref{contain0}.
\end{proof}

\begin{corollary}
Let $P$ be a Fano triangle with $ord(P)=\{a,b,c\}$ with $\\ a>b>c>0$.
Then, $P$ is of type $B_1$ if $a<b+c$.
\end{corollary}

\begin{remark}  
If the number of vertices of $P$ is equal to that of $B(P)$, then $$ord(B(P))=\Big\{\frac{g(i)}{I_p(v_{i-1}+v_{i})I_p(v_i+v_{i+1})}\Big\}.$$
\end{remark}

\section{Symmetric or  K\"ahler-Einstein Fano polytopes under B-transformation}

\subsection{Smooth case}\label{smoothcase}
We shall prove Theorem \ref{main-smooth} and Theorem \ref{BpreservesKE}.  

We deal with Fano polytopes for each fixed dimension. Note that there is only one Fano polytope $P = [-1, 1]$ of dimension one, whose associated toric variety is the projective line  $\mathbb{P}^1$, which is K\"ahler-Einstein. Note that $P$ is B-invariant, hence it is of type $B_\infty$.

There are $5$ smooth Fano polygons. (See \cite{N1} and \cite{KN}.)   It turns out that they are all of type $B_\infty$. More precisely, $\mathbb{P}^2$ is B-invariant and the other four smooth Fano polygons are $2$-periodic. Among them, the three Fano polygons corresponding to the projective plane $\mathbb{P}^2$, the Hirzebruch surface $\mathbb{F}_0$ of degree $0$ and the del Pezzo surface of degree $5$ are K\"ahler-Einstein.

Recall that smooth Fano polytopes  of dimension at most $6$ are completely classified in  \cite{GRDB} based on the algorithm in \cite{O}. From this,  we have a list of smooth Fano polytopes of   dimension at most $6$. By the help of computer (e.g., \cite{SageMath}), one can compute the number of smooth Fano polytopes of fixed dimension which is  of strict type $B_k$ for each given $k$.  In particular, there are $18$ smooth Fano $3$-polytopes  and we have the following table.  

\begin{center}
\begin{tabular}{|c||c|c|c|c|c|c||c|c|} 
\hline
\text{Strict type}     & $B_0$ & $B_1$ & $B_2$ & $B_3$ & $B_4$ & $B_{\infty}$ & \text{Total}& \text{KE} \\
\hline
\text{Number}    & 2 & 3& 7 &0 &1 & 5 & 18 & 5\\
 \hline
\end{tabular}
\end{center}

We emphasize  that the   $5$ Fano $3$-polytopes of type $B_\infty$ are precisely the $5$ K\"ahler-Einstein smooth Fano 3-polytopes, and they are  all $2$-periodic, hence it is of type $B_\infty$. Indeed, for a smooth Fano 3-polytope $P$ of type $B_5$, we have $B^2(P) = P$ except for one case (GRDB ID  = 18) in which case we have $B^3(P) = B(P)$.

Similarly, we compute the number of smooth Fano $4$-polytopes  of strict  type $B_k$ for each given $k$. Among $124$ smooth Fano $4$-polytopes, there are $14$  Fano $4$-polytope $P$ of type $B_4$,  $11$ of them being $2$-periodic, hence  of type $B_\infty$. The remaining three cases(GRDB ID: 53,54,55) are at least of type $B_{20}$. We expect that those three cases are of type $B_\infty$ but we do not know the proof.
$$\begin{array}{|c||c|c|c|c|c|c|c|c|c||c|c|}
\hline
\text{Strict type}     & B_0 &B_1 &B_2 &B_3&B_4&\ldots&B_{19}&B_{20}&B_\infty &\text{Total}&\text{KE} \\
\hline
\text{Number}    & 28 & 33 & 47  & 2 & 0&\ldots & 0 &\leq 3  & \geq 11 & 124 & 12\\
 \hline
\end{array}$$
There are precisely  $12$ smooth K\"ahler-Einstein Fano 4-polytopes where $11$ of them are of $B_\infty$ and the other one $P$ has GRDB ID 54. We do not know whether $P$ is periodic. 
But it seems that $P$ is pseudo-periodic, which is not the case for the other two Fano polytopes(GRDB ID:53,55) of type $B_{20}$. 
 
\begin{question}\label{pseudoperiodic}
Let $P$ be a smooth $4$-Fano polytope of type $B_\infty$. If $P$ is pseudo-periodic, then is $P$ K\"ahler-Einstein?
\end{question}

There are $866$ smooth Fano 5-polytopes. We can similarly compute the number of  smooth Fano $5$-polytopes of  strict type $B_k$ for each given $k$. It is worthwhile to note that the $23$ Fano polytopes of type $B_\infty$ are  precisely the $23$ K\"ahler-Einstein smooth Fano 5-polytopes.
$$\begin{array}{|c||c|c|c|c|c|c||c|c|}
\hline
\text{Strict type}     & B_0 &B_1 &B_2 &B_3&B_4&B_\infty &\text{Total}&\text{KE} \\
\hline
\text{Number}    & 342 & 278 & 215  & 4 & 4 &  23 &866 & 23 \\
 \hline
\end{array}$$


There are $7622$ smooth Fano 6-polytopes. In this case, there are $88$ smooth Fano 6-polytopes of type $B_3$, but we did not completely determine how many of them are of strict type. See Remark \ref{dim6}. The $23$ Fano $6$-polytopes of type $B_\infty$ are all $2$-periodic. In fact, each $P$ of them  satisfies $B^2(P)=B^4(P)$.
$$\begin{array}{|c||c|c|c|c|c||c|c|}
\hline
\text{Strict type}     & B_0 &B_1 &B_2 &B_3 & B_{\infty} & \text{Total} & \text{KE} \\
\hline
\text{Number}    & 3884 & 2510 & 1140  & \leq 65 & \geq 23 & 7622  &   51  \\
 \hline
\end{array}$$

\begin{remark}\label{dim6}
For a smooth Fano $6$-polytope $P$, $B^3(P)$ has a huge number of  vertices and facets. For example, the Fano $6$-polytope $P$ with GRDB ID 1787 has $12$ vertices and $80$ facets and  is of type $B_3$, but $B^3(P)$ has  $2772$ vertices and $1614$ facets. According to our estimation, it takes more than one month to precisely determine the number of smooth Fano $6$-polytopes of strict type $B_3$, so we have stopped the computation.
\end{remark}

Now, the following proposition will complete the proof of Theorem \ref{main-smooth}.
 
\begin{proposition}
$Q_1$, $Q_2$ and $Q_3$ in Theorem \ref{NP-theorem} are all of   type $B_\infty$.  
\end{proposition}

\begin{proof}
By computation, we see that $B^2(Q_1) = B^4(Q_1)$, $B^2(Q_2) = B^4(Q_2)$ and $B^2(Q_3) = B^4(Q_3)$, so $Q_1$, $Q_2$ and $Q_3$  are all periodic, hence they are of type $B_\infty$. 
\end{proof}

Since we have computed $B(P)$ for each  K\"ahler-Einstein Fano polytope $P$ of dimension at most $6$ or $Q_1$, $Q_2$ and $Q_3$ in Theorem \ref{NP-theorem}, it is easy to compute the barycenter of $B(P)$. This proves Theorem \ref{BpreservesKE}.

\subsection{Singular case} 
We shall prove Theorem \ref{main-singular}.

\subsubsection{Symmetric Fano polygons}

Consider the Fano polygon
$$S_{m,n}=conv\{(m+1,-m),(-m,m+1),(-n-1,n),(n,-n-1)\},$$
which is symmetric by \cite[Proposition 3.2]{HKim}.

\begin{lemma}\label{1-refl}
For every non-negative integers $m$ and $n$, the symmetric Fano polygon $S_{m,n}$ is of type $B_\infty$. Moreover, $B^k(S_{m,n})$ is  K\"ahler-Einstein for every $k \geq 1$.
\end{lemma}

\begin{proof} 
For any positive integer $t$, we have  $$B^{2t-1}(S_{m,n})=conv\{(-1,1),(-1,-1),(1,-1),(1,1)\}$$ and $$B^{2t}(S_{m,n})=conv\{(1,0),(0,1),(-1,0),(0,-1)\} = S_{0,0}.$$ 
Both of them are K\"ahler-Einstein, from which the result follows.
\end{proof}

\begin{remark} \label{Smnauto}
Note that $B(S_{m,n})$ always admits a rotation even though $S_{m,n}$ may not have one.  It is easy to see that $B^t(S_{m,n})$  admits a reflection for every positive integer $t$. 
\end{remark}

\begin{lemma}\label{I_p}
Let $P$ be a Fano polygon of type $B_1$. If $P$ admits a non-trivial automorphism $\sigma$, then 
 $I(v_i+v_{i+1})=I(\sigma(v_i)+\sigma(v_{i+1}))$ for all adjacent vertices $v_i$ and $v_{i+1}$ of $P$.
\end{lemma}

\begin{proof}
Since the lattice point $\frac{v_i+v_{i+1}}{I(v_i+v_{i+1})}$ is primitive, so is $$\sigma \Big(\frac{v_i+v_{i+1}}{I(v_i+v_{i+1})}\Big)=\frac{\sigma(v_i)+\sigma(v_{i+1})}{I(v_i+v_{i+1})}$$ by Lemma \cite[Lemma 2.2]{HKim}, from which the result follows. 
\end{proof}

\begin{proposition} \label{auto}
Let $P$ be a symmetric Fano polygon of type  $B_1$.
If $P$ admits a non-trivial rotation, then so does $B(P)$.
\end{proposition}

\begin{proof}
Let $v_1, \dots ,v_n $ be all vertices of $P$ written in counterclockwise order and $\sigma$ be a non-trivial rotation of  $P$.   
By Lemma \ref{I_p}, $I_p(v_i+v_{i+1})=I_p(\sigma(v_i)+\sigma(v_{i+1}))$ for all adjacent vertices $v_i$ and $v_{i+1}$ of $P$.
Let $$w_i=\frac{v_i+v_{i+1}}{I_p(v_i+v_{i+1})}$$ for every $i=1,\dots,n$ and  $W=\{w_1, \dots ,w_n\}$ be an ordered set. 
Then, $\sigma(W)=W$ and $B(P)=conv(W)$. Note that a point of $W$ need not be a vertex of $B(P)$ anymore in general. 
It is enough to show that $\sigma$ preserves the vertex set of $B(P)$. 
Again, it is enough to show that if $w_i$ is not a vertex of $B(P)$ then so is not $\sigma(w_i)$.

Suppose that $w_i$ is not a vertex of $B(P)$.
Let $w_j$ and $w_k$ be the vertices of $B(P)$ closest to $w_i$. 
Define
$$f(w_i)=det(w_{j},w_i)+det(w_{i},w_{k})+det(w_{k},w_{j})$$ where $j < i< k$. Then, $f(w_i) < 0$. Since $\sigma$ is a rotation,  $$f(\sigma(w_i))=det(\sigma(w_{j}),\sigma(w_i))+det(\sigma(w_{i}),\sigma(w_{k}))+det(\sigma(w_{k}),\sigma(w_{j}))<0.$$
Thus, by \cite[Proposition 5]{Su}, $\sigma(w_i)$ is not a vertex of $W$.   
\end{proof}

\begin{corollary}
If a matrix $\sigma$ induces an automorphism of a symmetric Fano polygon $P$,  then it also induces an automorphism of $B(P)$.
\end{corollary}

\begin{proof}
The result follows from Remark \ref{Smnauto} and Proposition \ref{auto} by \cite[Theorem 3.3]{HKim}.
\end{proof}

\subsubsection{K\"ahler-Einstein Fano polygons}
We first consider the triangle case.
 
\begin{theorem}\label{tri-B_inf}  
If $P$ is a K\"ahler-Einstein Fano triangle, then $P$ is of type $B_\infty$. Moreover, $B^{2s+1}(P)=-P$ and $B^{2s}(P)=P$ for all non-negative integer $s$. Hence, $B^s(P)$ is K\"ahler-Einstein   for every non-negative integer $s$.
\end{theorem}

\begin{proof} 
By Proposition \cite[Proposition 3.10]{HKim},  we may assume that $$P=conv\{(a,-b),(0,1),(-a,b-1)\}.$$
Then, it is easy to see that $B(P)=conv\{(a,-b+1),(-a,b),(0,-1)\}=-P.$ The rest follows immediately. 
\end{proof}

\begin{example}\label{inf-nonKE}
The converse of Theorem \ref{tri-B_inf} does not hold in general i.e., there exists a  Fano triangle of type $B_\infty$ that is not K\"ahler-Einstein. 
Let $$P=conv\{(-25,-12),(-5,-6),(25,14)\}.$$
Then,  it is easy to see that $P$ is not K\"ahler-Einstein. However, since  $$B(P)=conv\{(-5,-3),(5,2),(0,1)\}$$ is K\"ahler-Einstein, by Theorem \ref{tri-B_inf}, $B(P)$ is of type $B_\infty$, thus so is $P$.
\end{example}

We need the following technical lemma to prove Theorem \ref{polygon-KE}. 

\begin{lemma}\label{area} 
Let $P$ be a Fano triangle with vertices $(a,-b),(0,1),(-c,d)$
written in counterclockwise order where  $a,b,c,d$ are positive integers. 
Let $Q$ be the convex hull generated by $P$ and $v$ be a primitive lattice point   that lies below the line spanned by $(a,-b)$ and $(c,-d)$. 
Then,  we have the following.
\begin{enumerate}
    \item The dual triangle $P^\ast$ of $P$ has vertices $w_1,w_2 ,w_3$ written in counterclockwise order such that  $w_1$ and $w_2$ are strictly below the $x$-axis and $w_3$ is strictly above the $x$-axis.
    \item The dual polygon $Q^\ast$ of $Q$ has vertices $w_1,w_2,w'_3,w'_4$ written in counterclockwise order such that  $w'_3$ lies on the line segment connecting $w_2$ and $w_3$; and $w'_4$ lies on the line segment connecting $w_1$ and $w_3$.  
\end{enumerate}
\end{lemma}

\begin{proof}
Write  $v_1=(a,-b),v_2=(0,1)$, $v_3=(-c,d)$, and  $v=(x,y)$. 
Then, we can compute that 
$$P^\ast=conv\Big\{\Big(\frac{-b-1}{a},-1\Big),\Big(\frac{1-d}{c},-1\Big),\Big(\frac{b+d}{bc-ad},\frac{a+c}{bc-ad}\Big)\Big\}$$ and

$$Q^\ast=conv\Big\{ \Big(-\frac{b+1}{a},-1\Big),\Big(\frac{1-d}{c},-1\Big),\Big(\frac{y-d}{cy+dx},-\frac{c+x}{cy+dx}\Big),  \Big(-\frac{b+y}{bx+ay},\frac{x-a}{bx+ay}\Big)\Big\}.$$

Take $w_1=(\frac{-b-1}{a},-1)$, $w_2=(\frac{1-d}{c},-1)$,
$w_3=(\frac{b+d}{bc-ad},\frac{a+c}{bc-ad})$, $w'_3=(\frac{y-d}{cy+dx},-\frac{c+x}{cy+dx})$ and $w'_4=(-\frac{b+y}{bx+ay},\frac{x-a}{bx+ay})$.  Since $v_1$, $v_2$, $v_3$ are in counterclockwise order, we have $bc-ad > 0$. This proves (1). 

Since $Q$ is a Fano polygon, $$ord(v_3,v)+ord(v,v_1)+ord(v_1,v_3)=-cy-dx-bx-ay-(bc-ad)>0.$$
Thus, $$ord(w'_3,w_3)=\frac{-cy-dx-bx-ay-(bc-ad)}{(-cy-dx)(bc-ad)}>0$$ and $$ord(w_3,w'_4)=\frac{-cy-dx-bx-ay-(bc-ad)}{(bc-ad)(-bx-ay)}>0.$$
So, $w'_3,w_3,w'_4$ are written in counterclockwise order. 
Let $S(w_i,w_j)$ be the slope between $w_i$ and $w_j$.
Then, $$S(w_2,w_3)=S(w_2,w'_3)=\frac{c}{d} \ \ \ \ and \ \ \ \ S(w_1,w'_4)=S(w_1,w_3)=\frac{a}{b}.$$
Thus, we see that $w'_3$ lies on the line segment connecting $w_2$ and $w_3$,  and $w'_4$ lies on the line segment connecting $w_1$ and $w_3$.

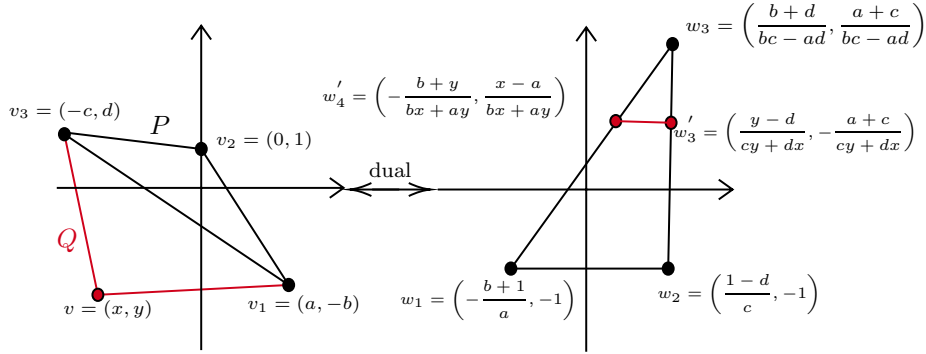
\tikzset{every picture/.style={line width=0.75pt}} 

\begin{tikzpicture}[x=0.75pt,y=0.75pt,yscale=-1,xscale=1]

\draw    (325.43,141.03) -- (327.64,27.73) ;
\draw    (246.11,141.03) -- (298.55,68.15) -- (327.64,27.73) ;
\draw [color={rgb, 255:red, 208; green, 2; blue, 27 }  ,draw opacity=1 ]   (134.02,149.3) -- (37.8,154.26) ;
\draw [color={rgb, 255:red, 208; green, 2; blue, 27 }  ,draw opacity=1 ]   (20.91,73.21) -- (32.28,127.74) -- (37.8,154.26) ;
\draw    (20.91,73.21) -- (134.02,149.3) ;
\draw  (17.06,100.3) -- (161.75,100.3)(89.91,18.42) -- (89.91,182.17) (154.75,95.3) -- (161.75,100.3) -- (154.75,105.3) (84.91,25.42) -- (89.91,18.42) -- (94.91,25.42)  ;
\draw  [fill={rgb, 255:red, 0; green, 0; blue, 0 }  ,fill opacity=1 ] (87.2,80.66) .. controls (87.2,78.94) and (88.43,77.55) .. (89.95,77.55) .. controls (91.47,77.55) and (92.71,78.94) .. (92.71,80.66) .. controls (92.71,82.37) and (91.47,83.76) .. (89.95,83.76) .. controls (88.43,83.76) and (87.2,82.37) .. (87.2,80.66) -- cycle ;
\draw  [fill={rgb, 255:red, 0; green, 0; blue, 0 }  ,fill opacity=1 ] (131.27,149.3) .. controls (131.27,147.59) and (132.5,146.2) .. (134.02,146.2) .. controls (135.54,146.2) and (136.78,147.59) .. (136.78,149.3) .. controls (136.78,151.01) and (135.54,152.4) .. (134.02,152.4) .. controls (132.5,152.4) and (131.27,151.01) .. (131.27,149.3) -- cycle ;
\draw  [fill={rgb, 255:red, 0; green, 0; blue, 0 }  ,fill opacity=1 ] (18.16,73.21) .. controls (18.16,71.5) and (19.39,70.11) .. (20.91,70.11) .. controls (22.43,70.11) and (23.67,71.5) .. (23.67,73.21) .. controls (23.67,74.93) and (22.43,76.31) .. (20.91,76.31) .. controls (19.39,76.31) and (18.16,74.93) .. (18.16,73.21) -- cycle ;
\draw  [fill={rgb, 255:red, 208; green, 2; blue, 27 }  ,fill opacity=1 ] (35.05,154.26) .. controls (35.05,152.55) and (36.28,151.16) .. (37.8,151.16) .. controls (39.33,151.16) and (40.56,152.55) .. (40.56,154.26) .. controls (40.56,155.97) and (39.33,157.36) .. (37.8,157.36) .. controls (36.28,157.36) and (35.05,155.97) .. (35.05,154.26) -- cycle ;
\draw    (21.65,72.39) -- (89.95,80.66) ;
\draw    (134.02,149.3) -- (89.95,80.66) ;
\draw  (209.2,101.12) -- (357.93,101.12)(284.09,19.25) -- (284.09,183) (350.93,96.12) -- (357.93,101.12) -- (350.93,106.12) (279.09,26.25) -- (284.09,19.25) -- (289.09,26.25)  ;
\draw  [fill={rgb, 255:red, 0; green, 0; blue, 0 }  ,fill opacity=1 ] (243.35,141.03) .. controls (243.35,139.32) and (244.59,137.93) .. (246.11,137.93) .. controls (247.63,137.93) and (248.86,139.32) .. (248.86,141.03) .. controls (248.86,142.74) and (247.63,144.13) .. (246.11,144.13) .. controls (244.59,144.13) and (243.35,142.74) .. (243.35,141.03) -- cycle ;
\draw  [fill={rgb, 255:red, 0; green, 0; blue, 0 }  ,fill opacity=1 ] (322.68,141.03) .. controls (322.68,139.32) and (323.91,137.93) .. (325.43,137.93) .. controls (326.95,137.93) and (328.19,139.32) .. (328.19,141.03) .. controls (328.19,142.74) and (326.95,144.13) .. (325.43,144.13) .. controls (323.91,144.13) and (322.68,142.74) .. (322.68,141.03) -- cycle ;
\draw  [fill={rgb, 255:red, 0; green, 0; blue, 0 }  ,fill opacity=1 ] (324.88,27.73) .. controls (324.88,26.01) and (326.11,24.63) .. (327.64,24.63) .. controls (329.16,24.63) and (330.39,26.01) .. (330.39,27.73) .. controls (330.39,29.44) and (329.16,30.83) .. (327.64,30.83) .. controls (326.11,30.83) and (324.88,29.44) .. (324.88,27.73) -- cycle ;
\draw  [fill={rgb, 255:red, 208; green, 2; blue, 27 }  ,fill opacity=1 ] (296.24,66.6) .. controls (296.24,64.88) and (297.47,63.5) .. (298.99,63.5) .. controls (300.51,63.5) and (301.75,64.88) .. (301.75,66.6) .. controls (301.75,68.31) and (300.51,69.7) .. (298.99,69.7) .. controls (297.47,69.7) and (296.24,68.31) .. (296.24,66.6) -- cycle ;
\draw  [fill={rgb, 255:red, 208; green, 2; blue, 27 }  ,fill opacity=1 ] (324.15,67.42) .. controls (324.15,65.71) and (325.38,64.32) .. (326.9,64.32) .. controls (328.42,64.32) and (329.66,65.71) .. (329.66,67.42) .. controls (329.66,69.14) and (328.42,70.53) .. (326.9,70.53) .. controls (325.38,70.53) and (324.15,69.14) .. (324.15,67.42) -- cycle ;
\draw    (246.11,141.03) -- (325.43,141.03) ;
\draw [color={rgb, 255:red, 208; green, 2; blue, 27 }  ,draw opacity=1 ]   (326.9,67.42) -- (298.99,66.6) ;
\draw    (176.15,101.12) -- (203.53,101.12) ;
\draw [shift={(205.53,101.12)}, rotate = 180] [color={rgb, 255:red, 0; green, 0; blue, 0 }  ][line width=0.75]    (10.93,-3.29) .. controls (6.95,-1.4) and (3.31,-0.3) .. (0,0) .. controls (3.31,0.3) and (6.95,1.4) .. (10.93,3.29)   ;
\draw    (189.78,101.12) -- (166.81,101.12) ;
\draw [shift={(164.81,101.12)}, rotate = 360] [color={rgb, 255:red, 0; green, 0; blue, 0 }  ][line width=0.75]    (10.93,-3.29) .. controls (6.95,-1.4) and (3.31,-0.3) .. (0,0) .. controls (3.31,0.3) and (6.95,1.4) .. (10.93,3.29)   ;

\draw (173.14,85.88) node [anchor=north west][inner sep=0.75pt]  [font=\footnotesize] [align=left] {dual};
\draw (62.27,62.49) node [anchor=north west][inner sep=0.75pt]   [align=left] {$P$};
\draw (15.53,119.5) node [anchor=north west][inner sep=0.75pt]   [align=left] {\textcolor[rgb]{0.82,0.01,0.11}{$Q$}};
\draw (111.81,153.12) node [anchor=north west][inner sep=0.75pt]  [font=\scriptsize] [align=left] {$\displaystyle v_{1} =( a,-b)$};
\draw (96.02,70.29) node [anchor=north west][inner sep=0.75pt]  [font=\scriptsize] [align=left] {$\displaystyle v_{2} =( 0,1)$};
\draw (-8.54,55.36) node [anchor=north west][inner sep=0.75pt]  [font=\scriptsize] [align=left] {$\displaystyle v_{3} =( -c,d)$};
\draw (19.04,155.65) node [anchor=north west][inner sep=0.75pt]  [font=\scriptsize] [align=left] {$\displaystyle v=( x,y)$};
\draw (188.74,145.42) node [anchor=north west][inner sep=0.75pt]  [font=\fontsize{0.63em}{0.76em}\selectfont] [align=left] {$\displaystyle w_{1} =\left(-\frac{b+1}{a} ,-1\right)$};
\draw (318.56,142.53) node [anchor=north west][inner sep=0.75pt]  [font=\fontsize{0.63em}{0.76em}\selectfont] [align=left] {$\displaystyle w_{2} =\left(\frac{1-d}{c} ,-1\right)$};
\draw (326.43,60.89) node [anchor=north west][inner sep=0.75pt]  [font=\fontsize{0.6em}{0.72em}\selectfont] [align=left] {$\displaystyle w^{'}_{3} =\left(\frac{y-d}{cy+dx} ,-\frac{a+c}{cy+dx}\right)$};
\draw (149.45,43.11) node [anchor=north west][inner sep=0.75pt]  [font=\fontsize{0.63em}{0.76em}\selectfont] [align=left] {$\displaystyle w^{'}_{4} =\left(-\frac{b+y}{bx+ay} ,\frac{x-a}{bx+ay}\right)$};
\draw (332.71,5.86) node [anchor=north west][inner sep=0.75pt]  [font=\scriptsize] [align=left] {$\displaystyle w_{3} =\left(\frac{b+d}{bc-ad} ,\frac{a+c}{bc-ad}\right)$};

\end{tikzpicture}

\end{proof}

It follows from Lemma \ref{area} that the $y$-coordinate of the barycenter of $P'$ is less than that of $P$ since $w_1$ and $w_2$ does not depend on the choice of $v$.
 
\begin{theorem}\label{polygon-KE}
 If $P$ is a K\"ahler-Einstein Fano polygon, then $P$ is of type  $B_1$. 
\end{theorem}

\begin{proof}
Let $P$ be a K\"ahler-Einstein Fano polygon. Assume that $P$ is not of type $B_1$. Then, by Proposition \ref{cond_B_1}, there exists a vertex $v_i$ of $P$ such that $g(i)\leq 0$. Up to a unimodular transformation, we may assume that $v_{i-1}=(a,-b), v_{i}=(0,1)$ and $v_{i+1}=(-c,d)$ for some positive integers $a,b,c,d$ where $v_{i-1}$, $v_i$,$v_{i+1}$ are adjacent vertices of $P$ in counterclockwise order.

Let $P'=conv\{v_{i-1},v_i,v_{i+1}\}$ and $(P')^\ast$ be its dual polygon. Note that the $y$-coordinate $y_{(P')^\ast}$ of the barycenter of $(P')^\ast$ is $$y_{(P')^\ast}= \frac{1}{3}\Big(   \frac{a+c}{bc-ad} -2\Big) =
\frac{(a+c-(bc-ad))-(bc-ad)}{3(bc-ad)}.$$
Since $g(i)=a+c-(bc-ad) \leq 0$, we have $y_{(P')^\ast} < 0$. 
By Lemma \ref{area}, if $P$ has another vertex $v$ other than $v_{i-1}$, $v_i$ and $v_{i+1}$, then $y_{P^\ast} \leq y_{(P')^\ast}<0$.
So, the barycenter of $P^\ast$ is not the origin.
Thus, $P$ is not K\"ahler-Einstein,  a contradiction.
\end{proof}

\begin{theorem}\label{polygon-symm}
Let $P$ be a Fano polygon. If $P$ is symmetric, then $P$ is of type $B_\infty$ and $B(P)$ is both symmetric and K\"ahler-Einstein. 
\end{theorem}
 
\begin{proof}
Let $P$ be a symmetric Fano polygon. If $P$ does not admit a rotation, then  $P = S_{m,n}$ with $m \neq n$ by  \cite[Theorem 3.3]{HKim}, thus the result follows from  Lemma \ref{1-refl}.
Now, we may assume that $P$ admits a rotation. Then, $P$ is K\"ahler-Einstein, and hence it is of type $B_1$ by Theorem \ref{polygon-KE}.  By Proposition \ref{auto},  $B(P)$ has an induced non-trivial rotation. Thus, $B(P)$ is also a symmetric Fano polygon of type  $B_1$. This argument shows that  $P$ is of type $B_\infty$. Moreover, $B^s(P)$ is K\"ahler-Einstein for every $s \geq 1$ by  \cite[Corollary 3.5]{HKim}.
\end{proof}

Now,  Theorem \ref{polygon-KE}, Theorem \ref{tri-B_inf}  and  Theorem \ref{polygon-symm}  prove Theorem \ref{main-singular}.

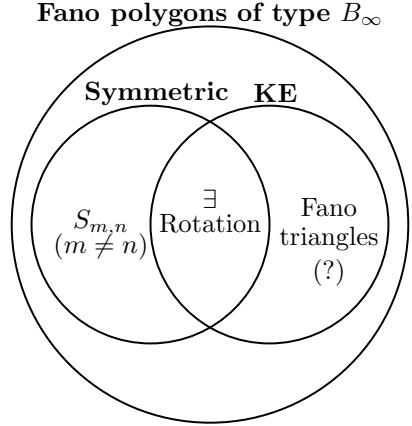
\begin{figure}[hb]\label{SymmKE-B}
\centering

\tikzset{every picture/.style={line width=0.75pt}} 
\begin{tikzpicture}[x=0.75pt,y=0.75pt,yscale=-1,xscale=1]
\begin{scope}
    \draw[red!30!white, draw = black] (-30,1) circle (60);
    \draw[blue!30!white, draw = black] (30,1) circle (60);
    \draw[blue!30!white, draw = black] (0,1) circle (100);
    \node at (0,-105) {\textbf{Fano polygons of type $B_\infty$}};
    \node at (-27,-65) {\textbf{Symmetric}};
    \node at (33,-65.5) {\textbf{KE}};
    \node at (-55,0) {$S_{m,n}$};
    \node at (-55,12) {($m \neq n$)};
    \node at (0,-12) {$\exists$};   
    \node at (0,0) {Rotation};
    \node at (60,-8) {Fano};
    \node at (60,8) {triangles};
    \node at (60,25) {(?)};
    \end{scope}
\end{tikzpicture}
    \caption{Symmetric/K\"ahler-Einstein polygons are of type $B_\infty$ assuming Conjecture \ref{triangle}.}
    \label{symm_KE_class}
\end{figure}

\begin{question}\label{BKE?}
Let $P$ be a K\"ahler-Einstein Fano polygon. Is $B(P)$ K\"ahler-Einstein?
\end{question}

If Question \ref{BKE?} is true, then Conjecture \ref{mainconj} is true for surfaces by Theorem \ref{main-singular}. Figure \ref{symm_KE_class} depicts Conjecture \ref{mainconj} for surfaces.

\subsubsection{Fano polygons of index at most $17$} 
We determine the number of Fano polygons of given strict type $B_k$ using the database \cite{GRDB} based on \cite{KKN} as we did in the smooth case. 

\begin{center}
\begin{tabular}{|c||c|c|c|c|c|c|c|c|c|c||c|c|}
\hline
\text{Index}     & $B_0$ &$B_1$ &$B_2$ &$B_3$&$B_4$&$B_5$&$B_6$&$B_7$&$B_8$&$B_\infty$ & \text{Total} & \text{KE} \\
\hline
1  & 3    & 0   & 0   & 1   & 0 & 0 & 0 &0&0& 12 &  16 & 5 \\
2  & 11   & 6   & 0   & 2   & 3 & 0 & 0 &0&0& 8 &  30 & 1 \\
3  & 35   & 20  & 10  & 2   & 4 & 2 & 0 &0&0& 26 &  99 & 4\\
4  & 35   & 25  & 10  & 2   & 3 & 0 & 0 &0&0& 16 &  91 & 2\\
5  & 100  & 63  & 35  & 10  & 6 & 2 & 0 &0&0& 34 &  250 & 7\\
6  & 126  & 93  & 49  & 25  & 8 & 8 & 1 &0&0& 69 &  379 & 4\\
7  & 186  & 100 & 73  & 9   & 7 & 1 & 0 &0&0& 53 &  429& 10\\
8  & 145  & 78  & 41  & 3   & 3 & 2 & 0 &0&0& 35 &  307 & 3\\
9  & 319  & 165 & 112 & 14  & 5 & 5 & 0 &0&0& 70 &  690  & 7\\
10 & 384  & 262 & 142 & 39  &15 & 3 & 1 &0&1& 69 &  916 & 4\\
11 & 472  & 190 & 155 & 28  & 6 & 0 & 0 &0&0& 88 &  939 & 12\\
12 & 535  & 325 & 185 & 74  & 26 & 9 & 2 &0&0& 123 &  1279 & 6\\
13 & 563  & 227 & 205 & 27  & 11 & 2 & 0 &0&0&  107 &  1142 & 17\\
14 & 725  & 423 & 261 & 34  & 20 & 2 & 1 &0&0&  79 &  1545 & 4\\
15 &1711  &1119 & 784 & 305 & 129& 24& 5 &1&0& 234 &  4312 & 14\\
16 & 564 & 237 & 137 &  9 & 9  & 5 & 0 &0&0& 69 &  1030 & 5\\
17 & 1007 & 353 & 330 &  53 & 9  & 1 & 0 &0&0& 139 &  1892&  19\\
 \hline
\end{tabular}
\end{center}

A direct computation based on the above proves Conjecture \ref{triangle} and Conjecture \ref{mainconj}  for Fano polygons of index at most $17$.
\begin{theorem}\label{index17}
Let $P$ be a   K\"ahler-Einstein Fano polygon  of index at most $17$. Then, 
    \begin{enumerate}
        \item $P$ is of type $B_\infty$ and $B(P)$ is again K\"ahler-Einstein.
        \item If $P$ is not symmetric, then $P$ is a triangle.
    \end{enumerate}
\end{theorem}

\section{Further discussions}
In this section, we discuss the number of Fano polytopes that are obtained by a sequence of B-transformations of a given Fano polytope, and study  the relations between the classes in Figure \ref{SymmKE-B} and the class of Fano polygons with barycenter zero.

\subsection{Orbits of a Fano polytope under B-transformation}
\begin{definition}
Let $P$ be a Fano polytope. Then, a Fano polytope $Q$ is called an \emph{orbit} of $P$ under B-transformation if it is obtained by applying a sequence of B-transformations of $P$, i.e., $Q = B^n(P)$ for some non-negative integer $n$.
\end{definition}
 
\begin{theorem}
Let $P$ be a symmetric Fano polygon that is not K\"ahler-Einstein or a K\"ahler-Einstein Fano triangle. Then, the number of orbits of $P$ under B-transformation  is  at most two.
\end{theorem}

\begin{proof}
It immediately follows from the proof of Lemma \ref{1-refl} and  Theorem \ref{tri-B_inf}.
\end{proof}

\begin{question}\label{orbitconj}
Let $P$ be a Fano polygon  of type $B_\infty$. Then, is the number of orbits of $P$ under B-transformation  at most two?
\end{question}

Conjecture \ref{triangle} gives a positive answer to Question \ref{orbitconj}.

\subsection{Fano polygons with barycenter zero}
We discuss the generalization of  the following proposition for K\"ahler-Einstein Fano triangles to arbitrary K\"ahler-Einstein Fano polygons.

\begin{proposition} \cite[Proposition 3.10]{HKim}
Let $P$ be a K\"ahler-Einstein Fano triangle. Then, $P$ has the origin as its barycenter.
\end{proposition}

\begin{theorem}\label{rot_bary0}
Let $P$ be a Fano polygon. If $P$ admits  a non-trivial rotation, then its barycenter  is the origin.  
\end{theorem}

\begin{proof}
Let $v$ be the barycenter of $P$. Since the set of vertices of $P$ is invariant under the action of rotations, it follows that $\sigma(v)=v$.
But, since the origin is the only fixed point of a rotation,
we conclude that $v$ is the origin.
\end{proof}

\begin{remark}\label{bary_S_mn}
Theorem \ref{rot_bary0} does not hold if $P$ admits no rotation.
Let $P=S_{m,n}$ for some non-negative integers $m$ and $n$. 
Then, the barycenter of  $P$ is $$\Big(\frac{m-n}{6(m+n+1)},\frac{m-n}{6(m+n+1)}\Big),$$ 
which is not the origin if $m \neq n$. Recall that $S_{m,n}$ is symmetric of type $B_\infty$ by Lemma \ref{1-refl}. 
\end{remark}

\begin{theorem}\label{KEbarycenter}
Assume that Conjecture \ref{triangle} holds. Then, every K\"ahler-Einstein Fano polygon has the origin as its barycenter. 
\end{theorem}

\begin{proof}
If $P$ has a non-trivial rotation, then $P$ is symmetric, so the barycenter of $P$ is the origin by Theorem \ref{rot_bary0}. If otherwise, $P$ is a K\"ahler-Einstein Fano triangle by Conjecture \ref{triangle}. Now the result immediately follows from   \cite[Proposition 3.10]{HKim}. 
\end{proof}

\begin{question}
Is it possible to remove the assumption in Theorem \ref{KEbarycenter}?  
\end{question}

Thanks to \cite{GRDB} based on \cite{KKN}, one can check that K\"ahler-Einstein Fano polygons are precisely the Fano polygons with barycenter zero if the index is at most 17. 
But there exist examples of Fano polygons with barycenter zero that are not K\"ahler-Einstein in higher index.

\begin{example} \label{bzero}
The following examples have the origin as their barycenters  but are not K\"ahler-Einstein. 
\begin{enumerate}
    \item Let $P_1=conv\{(-2, -1), (-1, 3), (1, 2), (2, -3)\}.$
Then, the barycenter of $P_1$ is the origin but $P_1$ is not K\"ahler-Einstein. 
Moreover, $P_1$ is of type $B_\infty$ and  the index of $P_1$ is $280$.  
\item Let $P_2=conv\{(-5, -4), (-5, 8), (5, 1), (8, -5)\}.$
Then, the barycenter of $P_2$ is the origin but $P_2$ is not K\"ahler-Einstein.
Moreover, $$B^5(P_2)=conv\{(-1, 0), (1, -1), (5, -3)\}, \text{ and }$$  $$B^6(P_2)=conv\{(0,-1),(3,-2),(4,-3)\}.$$
Thus, $P_2$ is of type $B_5$. Note that the index of $P_2$ is $270180$.   
\end{enumerate}
\end{example}

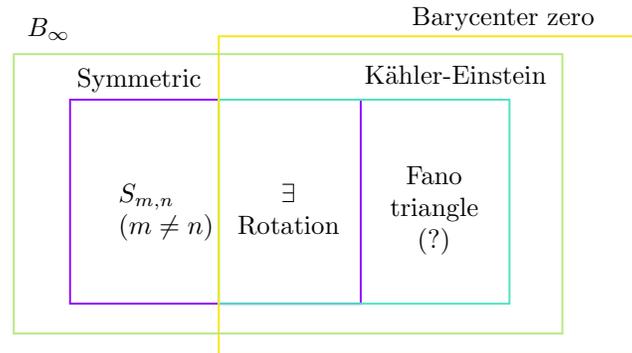
\begin{figure}[h]
    \centering
    
\tikzset{every picture/.style={line width=0.75pt}} 

\begin{tikzpicture}[x=0.75pt,y=0.75pt,yscale=-1,xscale=1]

\draw  [color={rgb, 255:red, 144; green, 19; blue, 254 }  ,draw opacity=1 ] (145,56.02) -- (291.61,56.02) -- (291.61,158.98) -- (145,158.98) -- cycle ;
\draw  [color={rgb, 255:red, 80; green, 227; blue, 194 }  ,draw opacity=1 ] (220,56.02) -- (366.61,56.02) -- (366.61,158.98) -- (220,158.98) -- cycle ;
\draw  [color={rgb, 255:red, 248; green, 231; blue, 28 }  ,draw opacity=1 ] (220,24.02) -- (434.61,24.02) -- (434.61,185.02) -- (220,185.02) -- cycle ;
\draw  [color={rgb, 255:red, 184; green, 233; blue, 134 }  ,draw opacity=1 ] (116.61,33.02) -- (393.29,33.02) -- (393.29,174.02) -- (116.61,174.02) -- cycle ;

\draw (147,39) node [anchor=north west][inner sep=0.75pt]   [align=left] {Symmetric};
\draw (316,8) node [anchor=north west][inner sep=0.75pt]   [align=left] {Barycenter zero};
\draw (292,37) node [anchor=north west][inner sep=0.75pt]   [align=left] {K\"ahler-Einstein};
\draw (168,97) node [anchor=north west][inner sep=0.75pt]   [align=left] {$\displaystyle S_{m,n}$ \\ $(m \neq n)$};
\draw (222,97) node [anchor=north west][inner sep=0.75pt]   [align=left] {\begin{minipage}[lt]{47.48406pt}\setlength\topsep{0pt}
\begin{center}
$\exists$ \\ Rotation
\end{center}

\end{minipage}};
\draw (303,88) node [anchor=north west][inner sep=0.75pt]   [align=left] {\begin{minipage}[lt]{36.163624000000006pt}\setlength\topsep{0pt}
\begin{center}
Fano\\triangle\\(?)
\end{center}

\end{minipage}};
\draw (122,13.4) node [anchor=north west][inner sep=0.75pt]    {$B_{\infty }$};
\end{tikzpicture}
    \caption{Hierarchy of  Fano polygons assuming Conjecture \ref{triangle}}
    \label{class}
\end{figure}

\begin{remark}
The Fano polygon $P$ in Example \ref{toKE} supports the example of a Fano polygon of type $B_\infty$ with nonzero barycenter that is neither symmetric nor K\"ahler-Einstein. The Fano polygon $P_1$ in Example \ref{bzero} (1) supports the example of a Fano polygon of type $B_\infty$ with barycenter zero that is neither symmetric nor K\"ahler-Einstein. The Fano polygon $P_2$ in Example \ref{bzero} (2) supports the example of a Fano polygon, not of type $B_\infty$,  with barycenter zero that is  neither symmetric nor K\"ahler-Einstein. 
\end{remark}

\bigskip {\bf Acknowledgements.}  
This research was supported by the Samsung Science and Technology Foundation under Project SSTF-BA1602-03.
 	 
\bigskip


\end{document}